# Automatic groups associated with word orders other than shortlex

Sarah Rees,
Department of Mathematics, University of Newcastle,
Newcastle NE1 7RU, UK,
e-mail: Sarah.Rees@ncl.ac.uk

**Abstract**

The existing algorithm to compute and verify the automata associated with an automatic group deals only with the subclass of shortlex automatic groups. This paper describes the extension of the algorithm to deal with automatic groups associated with other word orders (the algorithm has now been implemented) and reports on the use of the algorithm for specific examples; in particular a very natural automatic (or asynchonously automatic) structure for the Baumslag-Solitar and related classes of groups (closely related to one described for some of those groups by Epstein et al.) is found from a wreath product order over shortlex.

AMS subject classifications: 20F10, 20-04, 68Q40, secondary classification: 03D40.

## 1 Introduction

The definition of an automatic group grew out of work of Cannon ([3]) in an attempt to formalise properties observed in the fundamental groups of compact hyperbolic manifolds. Briefly (a formal definition will be given below), an automatic group is a group satisfying certain finiteness conditions which allow the group to be described by a set of finite state automata; one finite state automaton defines a language, or set of words in the generators, which contains at least one representative of each group element, the others describe multiplication between words in that language. The book [6] provides a standard reference.

*The author would like to thank both the Lehrstuhl D für Mathematik of the RWTH-Aachen, and the New York Group Theory Cooperative for their hospitality during the period of the research described in this paper, Derek Holt and David Epstein for helpful discussions, and the referee for valuable comments. The visit to Aachen was funded by the European Union HCM grant no. ERBCHRXCT930418

Using the automata associated with an automatic group, the 'word problem' can be solved, and the group elements efficiently enumerated. The growth function of the language for the defined group defined by the first automaton can be computed (this will be the growth function of the group if the language consists of unique geodesic representatives of the group elements). Some group theoretical properties (such as torsion) can sometimes be read from the automata. Thus it is interesting not merely to know that a group is automatic, but also to be able to construct explicitly the associated automata, and hence to compute information about the group. In [7] and [12], descriptions are given of two variations of a procedure which successfully construct these automata for groups which are automatic with respect to the SHORTLEX word order (see below for a definition of this subclass); the Warwick packages **automata** and its recent successor **kbmag** implement these ideas. Further a procedure originally described by Epstein, and now developed and implemented by P. Wakefield ([24]) constructs a further SHORT automatic structure from the SHORTLEX automatic structure of a word hyperbolic group (and so verifies hyperbolicity). However, although many topologically interesting groups are known to be SHORTLEX or SHORT automatic, some classes have only been shown to be automatic with respect to other word orders, while for other groups, geometrical considerations, or the existence of a finite rewriting system with respect to a word order other than SHORTLEX make it interesting to look for different kinds of automatic structures.

The purpose of this paper is to describe a generalisation of the original SHORTLEX procedure to construct the automata associated with a wider class of automatic groups, and hence verify automaticity for these groups. The procedure is designed in as general way as possible; it can be used for any reduction word order for which two particular ingredients can be defined. Most of the procedure is valid also for orders which are translation invariant but not actually reduction orders.

The procedure was initially implemented in $C++$ within the New York based **Magnus** system ([15]) for computation in group theory, specifically for the weighted orders WTLEX and WTSHORTLEX; in this implementation, it makes use of C. Sims' **rkbp** implementation of the Knuth-Bendix procedure ([23]), and parts of D. Holt's **kbmag** package for Knuth-Bendix and automatic groups. It has now been implemented within the GAP language ([22]), currently specifically for the WTLEX order and the wreath product order over SHORTLEX; in this implementation it also makes use of parts of the **kbmag** package. The code is available on request.

The paper is divided into six sections, of which the first is this introduction. Section 2 contains the definition of an automatic group. Section 3 outlines the tools and basic steps of the procedure to verify automaticity; this is modelled on the Warwick SHORTLEX procedure, but at some points needs to be more general. Section 4 describes the general construction of a word acceptor, this being the part of the procedure which needs substantial modification. The construction in section 4 is dependent on the particular order only to the extent that it requires the definition of the two ingredients already mentioned, that is (in the terms of this paper) a *history* function $f$ and a finite *sufficient* set $H$ of *difference histories* for that order; section 5 explains how $f$ and $H$ may be defined for each of the WTLEX order, the WTSHORTLEX order, and wreath product orders over SHORTLEX. Finally, section 6 looks at some examples of applications of the procedure for the above orders. In particular the structures obtained for the Baumslag-Solitar and related classes of groups seem to be interesting.

## 2 Definition of an automatic group

Throughout this paper $G$ will denote a finitely presented group defined by the presentation
$$G = \langle g_1, g_2, \ldots, g_m \mid x_1 = y_1, \ldots, x_k = y_k \rangle.$$
We shall use the term *generator* of $G$ to mean any element of the set $\{g_1, \ldots, g_m, g_1^{-1}, \ldots, g_m^{-1}\}$, and *word* in the generators of $G$ to mean any string of symbols from the set $\{g_1, \ldots, g_m, g_1^{-1}, \ldots, g_m^{-1}\}$.

For words $w, w'$ in the generators of $G$, $w = w'$ will mean that $w$ and $w'$ are identical as strings, while $w =_G w'$ will mean that $w$ and $w'$ represent the same group element of $G$. Futher, $l(w)$ will denote the string length of $w$ (for example, $l(g_1^5 g_2^{-3}) = 8$) and $l_G(w)$ the string length of a shortest word $u$ with $u =_G w$. For any integer $i$, $w(i)$ will denote the prefix of $w$ of length $i$ if $i < l(w)$, and will otherwise denote the whole of $w$. We shall use the symbol $\epsilon$ to denote the word of length 0, and $e$ to denote the identity element of the group $G$ (which $\epsilon$ represents).

For the purposes of this paper, a *finite state automaton* $M$ is a finite directed graph, whose edges are labelled by symbols from an associated set known as its *alphabet*. Out of each state there is at most one edge labelled with any given alphabet symbol. One vertex (or *state*) is known as the *start* state, some of the others are *accept* states. A string on the alphabet is accepted (and so is in the *language*, $L(M)$, of $M$) if it labels a path from the start state to an accept state, and is otherwise rejected. Since transitions out of some states on some symbols may be undefined, an automaton satisfying these conditions is often known as a *partial deterministic automaton*. A set of words is called a *regular set* if it is the language of a finite state automaton.

A *2-stringed automaton* is a finite state automaton whose alphabet is a set of ordered pairs of symbols from a *base alphabet*. The language can be considered as a set of pairs of words over the base alphabet; pairs of words are read by the automaton by reading at each stage a pair of symbols, one from each string. If one word is longer than the other, then the shorter word is padded out with a *padding symbol*, which we shall denote by $\_$, until the two words are equal in length.

**Definition 2.1 (Automatic group)** $G = \langle g_1, g_2, \ldots, g_m \rangle$ *is automatic if*

**(i)** *there is a finite state automaton $W$ with alphabet $\{g_1, \ldots, g_m, g_1^{-1}, \ldots, g_m^{-1}\}$ whose language $L(W)$ contains at least one representative of each element of $G$.*

**(ii)** *there is a constant $K$ such that, if $v, w \in L(W)$, and $l_G(v^{-1}w) \leq 1$, then, for all $i$, $l_G(v(i)^{-1}w(i)) \leq K$.*

Condition (ii) is known as the *fellow traveller property*. It can be replaced in the definition (and often is) by the following condition:-

**(ii)*** *for each element $g$ of the set $\{e, g_1, \ldots, g_m, g_1^{-1}, \ldots, g_m^{-1}\}$, there is a 2-stringed finite state automaton $M_g$, with base alphabet $\{g_1, \ldots, g_m, g_1^{-1}, \ldots, g_m^{-1}\}$ whose language is the set of pairs of words $(v, w)$ from $L(W)$ such that $vg =_G w$.*

For the purposes of computational verification of an automatic structure, (ii)* is sometimes more appropriate than the more geometrical condition (ii).

$W$ is called the *word acceptor*, and $M_g$ the *multiplier for $g$*. The word acceptor and multipliers together form an *automatic structure* for the group, and $L(W)$ the associated *language*. Using the multipliers, any word can be reduced to a representative in the language.

The above definition of automaticity appears to be dependent on the choice of generators for $G$. In fact it is not (see [6]). It is straightforward to deduce that $G$ is finitely presented.

The class of automatic groups includes finite groups, free groups, free abelian groups, various small-cancellation groups ([9]), Coxeter groups ([2]), Artin groups of finite and large type ([4, 5, 18, 13]), Braid groups ([6], chapter 9), word-hyperbolic groups and so, in particular, (see [3]), all co-compact discrete groups of isometries of $n$-dimensional hyperbolic space. However the class includes, for example, no virtually nilpotent group which is not virtually abelian.

Many of the languages associated with automatic structures are naturally associated with an easily described word order.

**Definition 2.2 (Automatic with respect to a word order)** *Let $\preceq$ be any partial order on words. A word $v$ is said to be* irreducible *with respect to $\preceq$ if there is no word $u$ with $u =_G v$ and $u \prec v$ (that is with $u \neq v$ and $u \preceq v$).*

*A group $G$ is said to be* automatic with respect to $\preceq$ *if there is an automatic structure for $G$ for which the language is the set of all irreducible words with respect to $\preceq$.*

Free abelian groups, Coxeter groups, and all word hyperbolic groups are automatic with respect to SHORTLEX, or 'SHORTLEX automatic', (for SHORTLEX, $u \prec v$ if $u$ is shorter than $v$, or if $u$ and $v$ have the same length but $u$ precedes $v$ lexicographically, in relation to a fixed order on the group generators and their inverses).

However for some groups, other languages are known to give automatic structures. By [6], Theorem 3.4.5, word hyperbolic groups are SHORT automatic, that is automatic with language the set of all geodesic words; Papasoglu ( [17]) showed that this property characterises word hyperbolic groups. The fundamental groups of geometrically finite hyperbolic manifolds were proved by Epstein ([6], Chapter 11) to be WTLEX automatic (given some weight function on words, $u \prec v$ in WTLEX if $u$ has lower weight than $v$ or if $u$ and $v$ have the same weight and $u$ precedes $v$ lexicographically). A new result by Urs Lang ([14]) now proves these groups also to be SHORTLEX automatic. Small cancellation groups satisfying the conditions $C(6) - P$, $C(4) - T(4) - P$ or $C(3) - T(6) - P$, which are not in general word hyperbolic, are proved in [9] to be automatic, with a language derived from a set of geodesic words. Braid groups and more generally Artin groups of finite type are proved to be automatic with a language associated with a 'left greedy' word order, and large type Artin groups (extra large, by [18], large by [13]) with respect to an order based on a combination of length, lexicographical order, and the natural order of associated binary strings. Finally Lee Mosher ([16]) has proved that the mapping class group of a surface of finite type is automatic with a language which is decribed geometrically in terms of triangulations on the surface.

Note that, whenever $G$ is automatic it is automatic with respect to the partial order $\preceq$ defined by

$$u \preceq v \iff v = u \quad \text{or} \quad v =_G u, u \in L(G), v \notin L(G)$$

However we shall be interested in the case where $\preceq$ is a reduction order, that is, a translation invariant linear order ($u \preceq v$ implies both $au \preceq av$ and $ub \preceq vb$ for all strings $a, b$) with no infinite descending chains $v_1 \succ v_2 \succ \ldots v_n \succ \ldots$.

## 3 The basic steps of the procedure to verify automaticity

Most of the theory and procedure of this section is simply a translation of the SHORTLEX automatic procedure of [7] and [12] into more general language. The main difficulty is in the construction of the word acceptor, and the associated reduction procedure, where, although the same basic principle is used, a great deal more care needs to be taken.

## 3.1 The word difference machine

An automaton known as a *word difference machine*, described below, is crucial to the procedure we shall describe for constructing and verifying an automatic structure.

Suppose that $G$ is automatic with language $L$. The fellow traveller property implies that the set $\mathcal{D}_L$ of all minimal representatives (under any specified reduction order) of elements $v(i)^{-1}u(i)$, with $v, u \in L$ and $l_G(v^{-1}u) \leq 1$, is finite. The elements of $\mathcal{D}_L$ are known as *word differences*.

We can define a 2-stringed automaton with the elements of $\mathcal{D}_L$ as states, in a natural way. We consider the empty word $\epsilon$ to be the start state, all states to be accept states, and define transitions on pairs of generators. We define a transition from the start state to itself on every pair $(g, g)$. In addition, for $d, d' \in \mathcal{D}$ we define a transition on the pair $(g, h)$ from $d$ to $d'$ if for distinct words $v, u \in L$, with $l_G(v^{-1}u) \leq 1$, and with $i$-th generators $g, h$,

$$d =_G v(i-1)^{-1}u(i-1), \ d' =_G v(i)^{-1}u(i).$$

There are no other transitions. Where $v$ (or $u$) has length less than $i$, its prefix $v(i)$ (or $u(i)$) of length $i$ is interpreted in the above as the whole string, and $g$ (or $h$) is interpreted as the padding symbol.

We call the automaton thus defined the *difference machine* for $L$ (and we use the notation $\mathcal{D}_L$ for it, as well as for the underlying set of word differences). The procedure we shall describe hinges on the construction of a machine which estimates $D_L$, and which we shall also call a difference machine.

**Definition 3.1 (Difference machine)** *We define a difference machine for $G$ to be a 2-stringed automaton $\mathcal{D}$, with base alphabet $X$, satisfying the following properties:-*

**(i)** *all of the states of $\mathcal{D}$ are accept states;*

**(ii)** *each of the states of $\mathcal{D}$ is labelled by a word in the generators of $G$, known as a word difference;*

**(iii)** *the start state is labelled by the trivial word;*

**(iv)** *if $d$, $d'$ are word differences, and if there is a transition from the state corresponding to $d$ to the state corresponding to $d'$ under the pair $(g, h)$, then $g^{-1}dh =_G d'$;*

**(v)** *for all $g \in X$, there is a transition from the start state to the start state on $(g, g)$.*

If $(w_1, w_2)$ is accepted by $\mathcal{D}$, $d_{w_1, w_2}$ is defined to be the label of the target of $(w_1, w_2)$; condition (iv) ensures that $d_{w_1, w_2}$ is a word representing $w_1^{-1}w_2$.

Clearly $\mathcal{D}_L$ as defined above is a difference machine according to this definition.

Condition (v) might be omitted. It makes the descriptions which follow in this paper more straightforward; but from an implementation point of view makes the difference machine bigger than it needs to be, and hence those transitions may not actually be stored, although the algorithms need to take account of them.

For notational convenience, we shall also define, for any word difference $d \in \mathcal{D}$, an automaton $\mathcal{D}^d$ which has the same states and transitions as $\mathcal{D}$ but only the single accept state $d$. Thus, for example, for any pair of words $(w_1, w_2)$ accepted by $\mathcal{D}^\epsilon$, $w_1 =_G w_2$.

**Definition 3.2 (Fellow travellers)** *If $(w_1, w_2)$ is a pair of words accepted by a difference machine $\mathcal{D}$, then we shall say that $w_1$ and $w_2$ $\mathcal{D}$-fellow travel.*

Note that condition (v) in the definition of a difference machine ensures that whenever $w_1$ and $w_2$ $\mathcal{D}$-fellow travel then so do $cw_1$ and $cw_2$ for any word $c$. However $w_1 c$ and $w_2 c$ might not.

The label $d_{w_1, w_2}$ of the target of the pair $(w_1, w_2)$ in $\mathcal{D}$, is a word representing the product $w_1^{-1} w_2$. In theory it could happen that $w_1$ and $w_2$ $\mathcal{D}$-fellow travel but $w_2$ and $w_1$ do not. In practice it will not happen in this paper. For in this paper, we shall often need our difference machines to satisfy the conditions defined below; of these, inverse closure ensures that $w_2$ and $w_1$ $\mathcal{D}$-fellow travel whenever $w_1$ and $w_2$ do.

**Definition 3.3 (Inverse closure, substring closure, minimal labelling)** *Let $\mathcal{D}$ be a difference machine.*

*We say that $\mathcal{D}$ is* inverse closed *if for each word difference $d$, a word in $G$ representing the inverse of $d$ is also a word difference, and whenever there is a transition from $d$ to $d'$ under $(g, h)$ there is a transition from $d'$ to $d$ under $(g^{-1}, h^{-1})$.*

*We say that $\mathcal{D}$ is* substring closed *if for each word difference $d$ both the maximal prefix $d_1$ of $d$ and the maximal suffix $d_2$ of $d$ are word differences, and where $d = d_1 g = h^{-1} d_2$ (for generators $g$ and $h$), there is a transition from $d_1$ to $d$ under $(\_, g)$ and from $d_2$ to $d$ under $(h, \_)$.*

*Finally we say that $\mathcal{D}$ is* minimally labelled *(with respect to an order $\preceq$) if the word $w$ which labels the state $d$ is the least of all words $w_1^{-1} w_2$ such that $(w_1, w_2)$ labels a path to $d$.*

Suppose from now on that $\preceq$ is translation invariant. Note that, in this case, any one of the above properties can be verified for a given difference machine $\mathcal{D}$ in finite time, and further, if any of them fail, a difference machine $\mathcal{D}'$ can be constructed in finite time, whose language contains the language of $\mathcal{D}$, and in which those properties hold.

To describe our procedure, we need the concepts of reducibility and irreducibility relative to a difference machine $\mathcal{D}$

**Definition 3.4 ($\mathcal{D}$-reducible, $\mathcal{D}$-irreducible)** *If $u$, $v$ are words such that $u =_G v$, $u \prec v$ and $u = au'b$, $v = av'b$ for some $a$, $b$ and some pair $(v', u')$ accepted by $\mathcal{D}^\epsilon$, we say that $v$ is $\mathcal{D}$-reducible, with direct $\mathcal{D}$-reduction $u$; $w$ is then called a $\mathcal{D}$-reduction of $v$ if there is a chain $u_0 = v, u_1, \ldots, u_k = w$, with $k > 0$ such that each $u_i$ is a direct $\mathcal{D}$-reduction of $u_{i-1}$. We say that $v$ is $\mathcal{D}$-irreducible if it is not $\mathcal{D}$-reducible.*

*Further we say that $v$ is* minimally $\mathcal{D}$-reducible *if it is $\mathcal{D}$-reducible, but every proper subword of it is $\mathcal{D}$-irreducible. We say that $v$ is directly $\mathcal{D}$-reducible if $v$ is $\mathcal{D}$-reducible, with some direct $\mathcal{D}$-reduction which is $\mathcal{D}$-irreducible.*

If $\mathcal{D} = \mathcal{D}_L$, for some automatic group with language $L$, then $L$ can be defined as the set of all $\mathcal{D}$-irreducible words. That all words in $L$ are $\mathcal{D}$-irreducible is clear. To see that $L$ contains all $\mathcal{D}$-irreducible words, notice that if $v$ is a word not in $L$, then $v = v_1 x v_2$ for words $v_1, v_2$ and a generator $x$ such that $v_1 \in L$, but $v_1 x \notin L$. Then where $u \in L$ is equal to $v_1 x$ in $G$, certainly $u \prec v_1 x$ and $u$ and $v_1 x$ $\mathcal{D}_L$-fellow travel. So $v$ is $\mathcal{D}_L$-reducible with direct $\mathcal{D}_L$-reduction $u v_2$.

If $\mathcal{D} = \mathcal{D}_L$, then every minimally $\mathcal{D}$-reducible word is directly $\mathcal{D}$-reducible; this follows from the above and the fellow traveller property for $L$. However this will not in general be true.

## 3.2 Defining the automatic structure in terms of the difference machine

Suppose that $\preceq$ is a translation invariant linear order. If $G$ is automatic with respect to $\preceq$, with language $L$, then the automatic structure is completely determined by the automaton $\mathcal{D}_L$ and $\preceq$.

The word acceptor can be defined as the finite state automaton which accepts all $\mathcal{D}_L$-irreducible words; the construction of this is described in section 4. The multiplier $M_g$ can be defined as the automaton which accepts $(v, w)$ if $v$ and $w$ are both accepted by $W$ and $(v, w)$ is accepted by $\mathcal{D}_L^d$, where $d = d_{\epsilon, g}$ (a word difference which represents $g$). Thus, a standard construction builds $M_g$ out of the automata $W$ and $\mathcal{D}_L^d$.

$\mathcal{D}_L$ can also be used to reduce a word to a representative in $L$. Given a word $w$, we trace pairs of words $(v, u)$ (for each substring $v$ of $w$) through $\mathcal{D}_L$, and look for a pair $(v, u)$ which is accepted by $\mathcal{D}_L^\epsilon$ and for which $u \prec v$. If such a pair exists, $v$ reduces to $u$, and hence $w$ also reduces. Iteration gives complete reduction.

## 3.3 The procedure as a whole

A construction to compute and verify an automatic structure with respect to the SHORTLEX order (and hence prove a group SHORTLEX automatic) is described in [7] and [12]. Here we describe the procedure in greater generality, in order to apply it to other orders $\preceq$.

The essential differences are in the details of the construction of the word acceptor, and in the procedure to reduce a word to one equal to it and accepted by the word acceptor. However we shall first describe the procedure as a whole, for greater clarity.

Briefly, a finite set of word differences is accumulated and an associated difference machine built, which is an estimate for the machine $\mathcal{D}_L$; the automata $W$ and $M_g$ described above are constructed from this. If $G$ is automatic with respect to $\preceq$ and if further the difference machine is indeed $\mathcal{D}_L$, then these will be the automata of an automatic structure, but otherwise one of a series of tests will fail. Some test failures indicate word differences which have not so far been found, and allow a modification of the difference machine, and then the reconstruction of the automata. If the group is automatic with respect to $\preceq$, then provided that the initial phase to accumulate word differences is run for 'long enough', after a finite number of loops the automatic structure will be found. (For some automatic groups the initial phase cannot find all word differences, and so some looping is essential.) However, if the group is not automatic with respect to $\preceq$, or even if it is but the initial phase was not run for long enough, the procedure may never terminate.

In fact the procedure described in [7] and [12] for SHORTLEX uses two different estimates for $\mathcal{D}_L$, at different stages; a smaller difference machine is used for the construction of the word acceptor than for the construction of the multipliers. We choose not to do that.

For the full procedure we require $\preceq$ to be a reduction order; otherwise the first step, in which the Knuth-Bendix ([10]) procedure is run, does not make sense. However much of the procedure depends only on $\preceq$ being translation invariant.

The basic steps of the full procedure are as follows:-

**Step 1:** Where $G = \langle g_j \mid x_k = y_k \rangle$ and $y_k \preceq x_k$, initialise a set $\mathcal{R}$ of rewriting rules to be
$$\{x_k \to y_k,\ g_j g_j^{-1} \to 1,\ g_j^{-1} g_j \to 1\}.$$

Run the Knuth-Bendix procedure on $\mathcal{R}$ for some time.

At any stage, given the rewriting system $\mathcal{R}$, define a related difference machine $\mathcal{D}$ as follows. The states of $\mathcal{D}$ are all reductions (by $\mathcal{R}$) of words of the form $x(i)^{-1}y(i)$ such that $x \to y$ or $y \to x$ is a rule in $\mathcal{R}$. For generators $g$, $h$, there is a transition on the pair $(g, h)$ from $d$ to $d'$ if $g^{-1}dh$ reduces in $\mathcal{R}$ to $d'$. Similarly there a transition on $(g, \_)$ or $(\_, h)$ from $d$ to $d'$ if $g^{-1}d$ or $dh$ reduces in $\mathcal{R}$ to $d'$.

If the Knuth-Bendix procedure terminates with a finite, confluent set of rules, or if otherwise $\mathcal{D}$ appear to stabilise, move on to step 2.

**Step 2:** Construct a word acceptor $W$, as a finite state automaton which accepts all $\mathcal{D}$-irreducible words, and does not accept both $w$ and $w'$ if $w'$ is a $\mathcal{D}$-reduction of $w$. Jump to step 4.

**Step 3:** If this step is entered, $\mathcal{D}$ has just been modified. If $W$ now accepts words $w, w'$ such that $w'$ is a $\mathcal{D}$-reduction of $w$, go back to step 2. Otherwise, continue to step 4.

**Step 4:** For each generator $g$, construct $M_g$ to accept a string $(v, w)$ if and only if $v \in L(W)$, $w \in L(W)$, $(v, w) \in L(\mathcal{D})$, and the target word difference of the pair $(v, w)$ is $d_{\epsilon,g}$. Define $M_e$ to accept $(v, v)$ if and only if $v$ is accepted by $W$.

**Step 5:** For each $M_g$, use standard automaton operations to check if the set

$$\{v : \exists w, (v, w) \in L(M_g)\}$$

is equal to $L(W)$. If this is true for all $g$, continue to step 6. If not, then for some $g$ a word $v$ can be identified in $L(W)$ for which no $(v, w)$ is accepted by $M_g$. Where $w$ is the reduction of $vg$ computed using the difference machine $\mathcal{D}$, compute and reduce the word differences associated with the equation $vg = w$, add them and the appropriate transitions to $\mathcal{D}$, and go back to step 3.

**Step 6:** For each of the defining relators $r = g_{i_1} g_{i_2} \ldots g_{i_k}$ of $G$, and also for each of the words $r = xx^{-1}$, where $x$ is a generator, compute an automaton $M_r$, which accepts $(v, w)$ if for some $v = v_0, v_1, v_2, \ldots, v_k = w$, $(v_{j-1}, v_j)$ is accepted by $M_{g_{i_j}}$. The group is proved automatic with the machines $W$ and all the $M_g$ defining its automatic structure if for each such $r$, $M_r$ accepts the same language as $M_e$.

If the Knuth-Bendix procedure does not terminate, then there is some room for interpretation as to when the machine $\mathcal{D}$ is considered to be stable. If $G$ is automatic, then at some stage during the Knuth-Bendix procedure, $\mathcal{D}$ must contain all the information necessary to build the word acceptor. However it might never contain enough infomation to build correct multipliers. For when $(v, w)$ is accepted by the correct automaton $M_g$, it can easily happen that the reduction rule $vg \to w$ would never be found by the Knuth Bendix procedure, but is only a consequence of a sequence of such rules.

If the Knuth-Bendix procedure terminates with a finite, confluent set, then a word acceptor may be built directly from that set, rather than from the difference machine; it is simply defined to reject a word if and only if it contains the left hand side of a rewrite rule as a substring. However, even in that case, the difference machine $\mathcal{D}$ computed from the rewrite system may not be large enough to build correct multipliers; indeed the group may not be automatic.

Step 5 checks to see if the language defined by the word acceptor $W$ satisfies the fellow traveller property. If step 5 completes successfully, then the automata $W, M_e, M_{g_1}, M_{g_1^{-1}}, \ldots, M_{g_m}, M_{g_m^{-1}}$ form an automatic structure for a monoid $H$ of which $G$ is a quotient.

Step 6 checks to see if the relations of $G$ hold in $H$; if they do then $G$ and $H$ are isomorphic, and so $G$ is automatic with the given automatic structure. If they do not it is because the Knuth-Bendix procedure was halted prematurely, some word differences were missed, and so $W$ accepts two different words which are equal in the group. If the Knuth-Bendix procedure in step 1 terminated with a confluent set, then step 6 is unnecessary, and may be omitted.

## 4 The general construction of the word acceptor

### 4.1 The basic idea

Suppose that $\preceq$ is a translation invariant linear order, and that $\mathcal{D}$ is an inverse closed and substring closed difference machine. (Such a difference machine can easily be constructed from a given one, if that machine does not already have the required properties.) Provided that $\preceq$ satisfies certain conditions (that is, admits the definition of the two essential ingredients, $f$ and $H$, as defined below), the construction we shall describe below terminates with a finite state automaton $W$ which accepts all $\mathcal{D}$-irreducible words for $\preceq$. $W$ may also accept some words which are $\mathcal{D}$-reducible, but will not accept both $w$ and $w'$ if $w'$ is a $\mathcal{D}$-reduction of $w$. If $L(\mathcal{D}) \supseteq L(\mathcal{D}_L)$, where $L$ is the language of an automatic structure, then $L(W)$ must be the set of all $\mathcal{D}$-irreducible words.

Where $\preceq$ is the SHORTLEX word order, it is straightforward to construct a machine $W$ which recognises precisely the $\mathcal{D}$-irreducible words. For in this case, there is a finite state automaton $GT$ which accepts a pair $(v, u)$ of words precisely when $v \succ u$. Then $W$ can be defined as the machine which rejects $w$ if for some prefix $v$ of $w$ there is a word $u$ such that $(v, u)$ is accepted both by $GT$ and by $\mathcal{D}^\epsilon$, more precisely

$$L(W) = (X^* \setminus \{v : \exists u, (v, u) \in L(GT) \cap L(\mathcal{D}^\epsilon)\})X^*$$

Since the class of finite state automata is closed under each of the logical operations *and*, *not* and $\exists$, and there are terminating algorithms to perform those operations, there is a clear terminating algorithm to compute the word acceptor in this case.

However, for a general word order $\preceq$ the set of pairs $(v, u)$ with $v \succ u$ is not a regular language, and so there is no finite state automaton to recognise it; instead we define a mechanism (equivalent to a finite state automaton, although we describe it in a different language) to recognise a large enough subset of that language for our purposes.

## 4.2 The ingredients of the general construction

There are two basic ingredients for the general construction, a *history* function $f$, and a set $H$ of *difference histories*, which is both finite and *sufficient* for reduction. All these terms will now be explained.

First we look for a way to encode the essential information for the comparison of two words under $\succ$. For this we use a *history function*.

**Definition 4.1 (History function)** *Let $\preceq$ be a translation invariant linear order on the set of words over some alphabet $A$. We define a* history function *for $\preceq$ to be a function $f$, defined on pairs $(v, u)$ of distinct words over $A$ with $l(v) \geq l(u)$, such that*

**(a)** *the sets of values of $f$ on $\{(v, u) : l(v) = l(u)\}$ and $\{(v, u) : l(v) > l(u)\}$ are disjoint;*

**(b)** *for any word $c$, and any $(v, u)$ with $l(v) \geq l(u)$, $f((v, u)) = f((cv, cu))$;*

**(c)** *for any $(v, u)$ with $v \neq u$ and $l(v) = l(u)$, knowledge of $f(v, u)$ is enough to decide*

  **(i)** *for any given words $a, b$ whether or not $va \succ ub$,*

  **(ii)** *for any given generators $g, h$, the value of $f(vg, uh)$;*

**(d)** *for any $(v, u)$ with $l(v) > l(u)$, knowledge of $f(v, u)$ is enough to decide*

  **(i)** *for any given word $a$, whether or not $va \succ u$,*

  **(ii)** *for any given generator $g$, the value of $f(vg, u)$.*

*We call the image of a pair of words or a set of pairs of words under $f$ its* history.

A little thought shows that a history function for SHORTLEX can be defined by

$$\begin{aligned} f(v,u) &= (1,0) \quad \text{if} \quad l(v) > l(u) \\ &= (0,1) \quad \text{if} \quad l(v) = l(u) \quad \text{and} \quad v \succ u \\ &= (0,0) \quad \text{if} \quad l(v) = l(u) \quad \text{and} \quad v \prec u \end{aligned}$$

In general if $\preceq$ is a word order for which the set of pairs $(v, u)$ with $v \succ u$ is regular, recognised by a finite state automaton $GT$, then the function $f$ which maps $(v, u)$ onto the target of $(v, u)$ in $GT$ is a history function.

For any $v$, in order to predict whether or not some word of the form $va$ is minimally $\mathcal{D}$-reducible, we need to consider pairs of words $(v, u)$, for which $u$ is no longer than $v$, and $\mathcal{D}$-fellow travels with $v$. If for some $a$, $va$ is minimally $\mathcal{D}$-reducible, then for some such $u$, some $ub$ is a direct $\mathcal{D}$-reduction of $va$; the pair $(d_{v,u}, f(v, u))$, carries enough information to recognise this. Hence, for any $\mathcal{D}$-reducible word $w$, we can locate a $\mathcal{D}$-reduction of $w$ by systematically computing all such pairs $(d_{v,u}, f(v, u))$ for which $v$ is a prefix of $w$, and for each such, where $g$ is the generator in $w$ following the prefix $v$, testing to see whether $vg$ reduces to some $ub$. An iteration of this procedure reduces $w$ to a $\mathcal{D}$-irreducible word $w'$, with $w' =_G w$.

This motivates our definition of a *difference history*.

**Definition 4.2 (Difference history)** *Let $\mathcal{D}$ be a difference machine, $\preceq$ a translation invariant word order, and $f$ a history function for $\preceq$. For any word $v$, we define a difference history $\delta$ for $v$ (with respect to $(\mathcal{D}, f)$) to be a pair $(d, f(v, u))$ such that $v \neq u$, $l(v) \geq l(u)$, $(v, u) \in L(\mathcal{D})$, and $d = d_{v,u}$. We define $f^{-1}(\delta)$ to be the set of all pairs $(v', u')$ of distinct words with $l(v') \geq l(u')$ associated with the word difference $d$ for which $f(v', u') = f(v, u)$. Where $H$ is a set of difference histories, we define $f^{-1}(H)$ analogously.*

**Definition 4.3 (Target of a difference history)** *Let $\delta = (d, f(v, u))$ be a difference history, and suppose that $g, h$ are generators.*

*If the target of $d$ under $(g, h)$ is a non-identity word difference $d'$, and if $l(v) = l(u)$, then the target of $\delta$ under $(g, h)$ is defined to be the difference history $(d', f(vg, uh))$, otherwise the target of $\delta$ under $(g, h)$ is undefined.*

*If the target of $d$ under $(g, \_)$ is a non-identity word difference $d''$, then the target of $\delta$ under $(g, \_)$ is defined to be the difference history $(d'', f(vg, u))$, otherwise the target of $\delta$ under $(g, \_)$ is undefined.*

When $f$ has finite range, we have basically defined the states and transitions of a finite state automata, whose states are the difference histories. However, in general we cannot find $f$ with finite range. If we are only interested in accumulating the tree of difference histories necessary to $\mathcal{D}$-reduce a single word, this is not a problem; however we need to be more careful in order to use the same mechanism to build a word acceptor, since otherwise the construction might not terminate.

We deal with this problem by restricting consideration to a proper subset of the set of difference histories; we hope to find an appropriate subset which is also finite.

**Definition 4.4** *Let $H$ be a set of difference histories for $f$.*

*Where $v$ is minimally $\mathcal{D}$-reducible, we say that $H$ carries a reduction for $v$ (to $u$) if $u$ is a $\mathcal{D}$-reduction of $v$ such that for every proper prefix $v(i)$ of $v$, $(v(i), u(i))$ is in $f^{-1}(H)$. More generally, where $v$ is any $\mathcal{D}$-reducible word, we say that $H$ carries a reduction of $v$ if it carries a reduction for some minimally $\mathcal{D}$-reducible subword of $v$.*

*$H$ is said to be* sufficient *for reduction if for any minimally $\mathcal{D}$-reducible word $v$*

**either** *$H$ carries a reduction for $v$*

**or** *$v$ is not directly $\mathcal{D}$-reducible, and $H$ carries a reduction for every direct $\mathcal{D}$-reduction of $v$.*

Of course the set of all difference histories is always sufficient, but we are in particular interested in sets which are both sufficient and finite.

### 4.3 Using $f$ and $H$ to build the word acceptor

Suppose from now on that we have a history function $f$, and a set of difference histories $H$ which is both sufficient and finite.

Briefly, we construct a word acceptor by building outwards from the start state. First the target of the start state under each generator $g$ is defined, either to be failure, or to be a set of difference histories in $H$, and those targets which are not the failure state are put on a queue to be dealt with. From then on, states are taken off the queue one at a time, and for each, its target under each generator is computed (as described below), as a set of difference histories; this may then be identified either as the set corresponding to an existing state, or as a new state. Each time a new state is found it is put onto the end of the queue. The word acceptor is complete when the queue is empty.

The procedure is designed so that the state $s$ of the word acceptor reached by tracing out a non-trivial accepted word $w$ corresponds to a set of difference histories for $w$. For any $g$, when the target of $s$ under $g$ is an accept state, it will correspond to a set of difference histories which are the targets of the difference histories associated with $s$.

We shall prove that the construction terminates, given that $H$ is finite, and yields a correct word acceptor, given that $H$ is sufficient.

For any word order $\preceq$ for which the language of pairs $(v, u)$ with $v \succ u$ is regular, the language of a finite state automaton $GT$, the construction below is a natural construction of the machine to recognise the language

$$(X^* \setminus \{v : \exists u, (v, u) \in L(GT) \cap L(\mathcal{D}^\epsilon)\})X^*,$$

when $f$ simply maps $(v, u)$ onto the target state of $(v, u)$ in $GT$, and $H$ is simply the set of all difference histories.

First we give the details of the construction.

**Where $s_0$ is the start state, for each generator $g$ :**

> define the target of $s_0$ under $g$ to be failure if for some (not necessarily non-trivial) word $z$ with $z \prec g$, $(g, z)$ is accepted by $\mathcal{D}^\epsilon$;
>
> otherwise define the target $t$ of $d$ under $g$ to correspond to the set of difference histories in $H$
>
> $$\{(d, f(g, h)) : d = d_{g,h} \in \mathcal{D}, h \quad \text{is a generator or} \quad \epsilon\}$$
>
> and add $t$ to the end of the queue.

**For each state $s$ on the queue in turn, for each generator $g$ :**

> define the target of $s$ under $g$ to be failure if either the target of the start state under $g$ is failure (in this case $g$ itself reduces), or for some difference history $\delta$ associated with $s$, with associated word difference $d$, one of the following situations occurs:-
>
> **either (a)** the target of $d$ under $(g, \_)$ is accepted by $\mathcal{D}^\epsilon$ (and so, $g^{-1}d =_G e$) and examination of $\delta$ shows that if $(v, u) \in f^{-1}(\delta)$, then $u \prec vg$,
>
> **or (b)** there exists a generator $h$ such that the target of $d$ under $(g, h)$ is accepted by $\mathcal{D}^\epsilon$ (and so, $g^{-1}dh =_G e$) and examination of $\delta$ shows that if $(v, u) \in f^{-1}(\delta)$, then $uh \prec vg$,
>
> **or (c)** there exists a generator $h$ such that the target of $d$ under $(g, h)$ is a word difference $d'$, and examination of $\delta$ shows that, where $d''$ is the word difference inverse to $d'$, if $(v, u) \in f^{-1}(\delta)$, then $uhd'' \prec vg$;
>
> otherwise define the target $t$ of $s$ under $g$ to correspond to the set of difference histories in $H$ which are targets under $(g, h)$, with $h$ a generator or $\epsilon$, of difference histories in the set corresponding to the state $s$, together with all difference histories in $H$ which are associated with the target in the word acceptor of $s_0$ under $g$.

**Proposition 4.5** *Let $\preceq$ be a translation invariant linear order, and let $\mathcal{D}$ be an inverse closed and substring closed difference machine. When $f$ is a history function for $(\preceq, \mathcal{D})$, and $H$ a finite sufficient set of difference histories for $f$, the above construction terminates and constructs an automaton which accepts all $\mathcal{D}$-irreducible words, and cannot accept both $w$ and $w'$ if $w'$ is the $\mathcal{D}$-reduction of $w$.*

PROOF: Termination is guaranteed by the finiteness of $H$; $H$ can thus only have finitely many subsets to correspond to states.

The construction ensures that the state $s$ which is the target of an accepted string $w$ is associated with a set of difference histories $\delta$ each of the form $(d, f(v, u))$ for some suffix $v$ of $w$. If $wg$ is rejected, then for some difference history $(d, (v, u))$ associated with $s$, and for $u'$ defined to be either $u$, $uh$ for some generator $h$, or $uhd''$, we have $u' =_G vg$ and $u' \prec vg$. Since $v$ is a suffix of $w$, we have $wg = avg$ for some $a$. Then $vg$ and $u'$ $\mathcal{D}$-fellow travel, since $v$ and $u$ do (where $u' = uhd''$ we use the fact that $\mathcal{D}$ is inverse and substring closed). Hence so do $wg$ and $au'$. So $wg$ is $\mathcal{D}$-reducible, with $\mathcal{D}$-reduction $au'$.

This proves that every string rejected by the automaton is $\mathcal{D}$-reducible, and hence that any $\mathcal{D}$-irreducible word must be accepted.

Now suppose that $w$ is any word for which $H$ carries a reduction. We show that $W$ rejects $w$.

The word $w$ can be written in the form $avb$ for some minimally $\mathcal{D}$-reducible $v$ for which $H$ carries a reduction. Then, there is some $u \prec v$ such that for each proper prefix $v(i)$ of $v$, the pair $(v(i), u(i))$ is in $f^{-1}(H)$. Let $v = g_{j_1} g_{j_2} \ldots g_{j_r}$, and let $u = h_{j_1} h_{j_2} \ldots h_{j_{r'}}$. Since $v$ is minimally $\mathcal{D}$-reducible, $g_{j_1} \neq h_{j_1}$. If $r' < r$, now define $h_{j_i}$ to be the padding symbol for $i > r'$, for ease of notation.

Let $s$ be the target of the start state under $a$, and for each $i \leq l(v)$, let $s_i$ be the target of the start state under $av(i)$. (We may assume without loss of generality that these are all defined, since otherwise certainly $w$ is rejected.)

Since each $(v(i), u(i))$ is in $f^{-1}(H)$, difference histories $\delta_i = (d_i, f(v(i), u(i))) \in H$, can be defined for each $i$, and $\delta_i$ is the target of $\delta_{i-1}$ under $(g_{j_i}, h_{j_i})$ for any $i > 1$. $\delta_1$ is associated with the target of the start state under $(g_{j_1}, h_{j_1})$ and hence also with $s_1$, as the target of the state $s$ under $(g_{j_1}, h_{j_1})$. Hence $\delta_i$ is associated with $s_i$ for each $i$, and in particular $\delta_{r-1}$ is associated with $s_{r-1}$.

The target of $s_{r-1}$ under $g_{j_r}$ is now defined to be the failure state by rule (a) if $r > r'$, by rule (b) if $r = r'$ and by rule (c) if $r < r'$.

Hence $av$ and also $w = avb$ is rejected by the automaton.

Now finally, suppose that $w$ is a $\mathcal{D}$-reducible word which is accepted by $W$, and suppose that $w'$ is a $\mathcal{D}$-reduction of $w$. Since $w$ is accepted, $H$ cannot carry a reduction for $w$. So, since $w'$ is a $\mathcal{D}$-reduction of $w$, and $H$ is sufficient, $H$ must carry a reduction for $w'$. Hence $w'$ is rejected by $W$. □

**Remarks**

1. The procedure has been described in its most basic form. For an efficient implementation, improvements are necessary. A substantial improvement can be obtained if we can define a partial order of difference histories, for which, if $\delta' = (d_{v,u'}, f(v, u'))$ precedes $\delta = (d_{v,u}, f(v, u))$, then whenever $va$ is minimally $\mathcal{D}$-reducible and $H$ carries a reduction of $va$ to $ub$, $H$ also carries a reduction of $va$ to some $u'b'$. Given such a partial order, the procedure can be modified so that a difference history $\delta$ is discarded (that is, not included in the set of such associated with a state of the word acceptor) if there is another difference history $\delta'$, preceding it in the order, which is also associated with that state.

2. If $\mathcal{D}$ were not substring closed, or $H$ were not sufficient, then the procedure would certainly terminate, but with a word acceptor which might accept a pair of $\mathcal{D}$-fellow travellers $w, w'$ representing the same group element. If this were the case, it would be recognised during the construction of the multipliers (that is, during step 4 of the full procedure). Hence, in some cases, in order to restrict the size of the computation, it might be appropriate to work with $\mathcal{D}$ and $H$ for which the conditions of substring closure and sufficiency are not clearly satisfied.

## 5  $f$ and $H$ for various word orders

In this section, appropriate $f$ and $H$ are described for the WTLEX, and WT-SHORTLEX orders and the wreath product order over SHORTLEX. In each case the sufficiency of $H$ is proved. The results of implementations of the procedure for each of these orders are discussed in section 6.

The following technical lemma will be useful for proving the sufficiency of some sets of difference histories. We state it in a rather general form so that we can apply it to each of the word orders considered in this paper.

**Lemma 5.1** *Suppose that $\preceq$ is a word order, and that $\mu$ is an additive positive valued function on words which is compatible with $\preceq$ (that is, such that $\mu(u) \leq \mu(v)$ whenever $u \preceq v$, and $\mu(v_1 v_2) = \mu(v_1) + \mu(v_2)$). Suppose that $\mathcal{D}$ is an inverse closed and substring closed difference machine. Let $w_1, w_2$ be words which $\mathcal{D}$-fellow travel, and let $d = d_{w_1,w_2}$ be the target word difference of $(w_1, w_2)$. Then if $w_2$ is $\mathcal{D}$-irreducible, and no longer than $w_1$,*

$$\mu(w_2) - \mu(w_1) \leq \mu(d)$$

PROOF: Since $w_2$ is no longer than $w_1$ and $\mathcal{D}$ is substring closed, $w_1 d$ and $w_2$ fellow travel. If $\mu(w_2) - \mu(w_1) > \mu(d)$, then

$$\mu(w_1 d) = \mu(w_1) + \mu(d) < \mu(w_2),$$

and so $w_1 d$ is a $\mathcal{D}$-reduction of $w_2$, contradicting the fact that $w_2$ is $\mathcal{D}$-irreducible. $\square$

### 5.1  The WTLEX word order

Now, suppose that $\preceq$ is the WTLEX word order. The order is defined by the rule that $u \prec v$ if $wt(u) < wt(v)$, or if $wt(u) = wt(v)$ and $u$ precedes $v$ lexicographically.

We can define a history function $f$, on pairs $(v, u)$ of distinct words with $l(v) \geq l(u)$, as follows:-
$$f(v, u) = (len, lex, wtd)$$

where

$len\ =\ 1$ if $l(v) > l(u)$ and $len = 0$ if $l(v) = l(u)$;
$lex\ =\ 1$ if $u$ precedes $v$ lexicographically, otherwise $lex = -1$;
$wtd\ =\ wt(v) - wt(u)$ if $l(v) = l(u)$, otherwise $wtd = min\{wt(v) - wt(u), 1\}$.

It is left to the reader to verify that $f$ has the appropriate properties.

**Proposition 5.2** *Let $H$ be the set of all difference histories $(d, (len, lex, wtd))$ such that*
$$-wt(d) \leq wtd \leq max_{d' \in \mathcal{D}} wt(d')$$

*Then $H$ is sufficient.*

PROOF: Suppose that $v$ is minimally $\mathcal{D}$-reducible, and suppose that $u$ is any reduction of $v$. Let $v(i)$ is any proper prefix of $v$. Then $l(v(i)) = i \geq l(u(i))$. If $l(u(i)) = i$, then we can apply lemma 5.1 with $w_1 = u(i)$ and $w_2 = v(i)$ and the word difference $d_{u(i),v(i)}$ to get

$$wt(v(i)) - wt(u(i)) \leq wt(d_{u(i),v(i)}) \leq max_{d' \in \mathcal{D}} wt(d')$$

If $l(u(i)) < i$, then the last component of $f(v(i), u(i))$ is defined to be at most 1.

Now first suppose that $v$ is directly $\mathcal{D}$-reducible. In that case let $u$ be a $\mathcal{D}$-reduction of $v$ which is $\mathcal{D}$-irreducible. Applying lemma 5.1 with $w_1 = v(i)$, $w_2 = u(i)$ and the word difference $d_{v(i),u(i)}$ we see that

$$wt(v(i)) - wt(u(i)) \geq -wt(d_{v(i),u(i)}).$$

So $(d_{v(i),u(i)}, f(v(i), u(i)))$ is a difference history in $H$. Hence in this case $H$ carries a reduction for $v$.

Now suppose that $H$ does not carry a reduction for $v$. In this case the above argument shows that $v$ cannot be directly $\mathcal{D}$-reducible. Now suppose that $u$ is a direct $\mathcal{D}$-reduction for $v$, and let $v'$ be the shortest prefix of $v$ for which, where $i = l(v')$, $u' = u(i)$ and $d = d_{v',u'}$, the difference history $(d, f(v', u'))$ is not in $H$. Then $wt(v') - wt(u') < -wt(d)$, so $wt(v'd) < wt(u')$, and thus $v'd \prec u'$. Since $\mathcal{D}$ is substring closed, $v'd$ and $u'$ must $\mathcal{D}$-fellow travel, so $v'd$ is a $\mathcal{D}$-reduction of $u'$, which is no shorter than $u'$. For any proper prefix $u'' = u'(j)$ of $u'$, $v'' = v'(j)$ is a proper prefix of $v'$, and so is certainly irreducible, so the difference history $(d_{u'',v''}, f(u'', v''))$ must be in $H$ by lemma 5.1. Hence $H$ carries a reduction for $u$. □

## 5.2 The WTSHORTLEX word order

Now, suppose that $\preceq$ is the WTSHORTLEX word order. The order is defined by the rule that $u \prec v$ if $wt(u) < wt(v)$, or if $wt(u) = wt(v)$ and $u$ precedes $v$ in SHORTLEX

We can define a history function $f$ almost exactly the same as for WTLEX

$$f(v, u) = (len, lex, wtd)$$

$$\begin{aligned}
len &= 1 \quad \text{if} \quad l(v) > l(u) \quad \text{and} \quad len = 0 \quad \text{if} \quad l(v) = l(u); \\
lex &= 1 \quad \text{if} \quad l(v) = l(u) \quad \text{and} \quad u \quad \text{precedes} \quad v \quad \text{lexicographically}, \\
lex &= -1 \quad \text{if} \quad l(v) = l(u) \quad \text{and} \quad v \quad \text{precedes} \quad u \quad \text{lexicographically, otherwise} \quad lex = 0 \\
wtd &= wt(v) - wt(u) \quad \text{if} \quad l(v) = l(u), \quad \text{otherwise} \quad wtd = min\{wt(v) - wt(u), 1\}.
\end{aligned}$$

We define a sufficient set $H$ just as for WTLEX.

## 5.3 Wreath product orders

For a wreath product order, each generator or inverse of a generator is assigned a level. The level, $\lambda(v)$ of a word $v$ is defined to be the maximum of the levels of all generators in it. Then, for any $j$, $\text{pfx}_n(v)$ is defined to be the maximal prefix of $v$ consisting only of symbols of weight at most $j$, and $\pi_j(v)$ is defined to be the projection of $\text{pfx}_j(v)$ onto its weight $j$ symbols (that is the word formed by deleting from $\text{pfx}_j(v)$ all symbols of weight other than $j$).

For a pair of words $v, u$ with no common prefix, the wreath product $\preceq$ over an order ORD is then defined (recursively) by the following rule:-

$$\begin{aligned}
u \prec v \iff \quad &\text{either} \quad \lambda(u) < \lambda(v), \\
&\text{or} \quad \lambda(u) = \lambda(v) = J, \quad \text{and} \quad \pi_J(u) <_{\text{ORD}} \pi_J(v), \\
&\text{or} \quad \lambda(u) = \lambda(v) = J, \quad \pi_J(u) = \pi_J(v), \quad \text{and} \quad \text{pfx}_{J-1}(u) \prec \text{pfx}_{J-1}(v).
\end{aligned}$$

Further, for any word $a$, $au \prec av$ if $u \prec v$.

Suppose that $\preceq$ is a wreath product order, over ORD, with levels from $1 \ldots N$. Suppose that $g$ is a history function for ORD. Then a little thought should show that the following defines a history function $f$ for $\preceq$.

Let $J = max\{\lambda(v), \lambda(u)\}$. Then, for distinct $v, u$ with $l(v) \geq l(u)$ and no common prefix, define

$$f(v, u) = (len, \lambda(v), \lambda(u), c_1, \ldots, c_{J-1}, c_J),$$

where $\quad len = 1 \quad \text{if} \quad l(v) > l(u), \quad = 0 \quad \text{if} \quad l(v) = l(u),$

and for each $j \leq J$,

$$\begin{align}
c_j &= 1 \quad \text{if} \quad \pi_j(u) <_{\text{ORD}} \pi_j(v) \quad \text{and} \quad (l(u) < l(v) \quad \text{or} \quad \lambda(u) > j); \\
c_j &= -1 \quad \text{if} \quad \pi_j(v) <_{\text{ORD}} \pi_j(u) \quad \text{and} \quad \lambda(v) > j; \\
c_j &= 0 \quad \text{if} \quad \pi_j(v) = \pi_j(u);
\end{align}$$

otherwise $c_j = (g(v_{j,1}, u_{j,1}), v_{j,2}, u_{j,2})$,

where $\pi_j(v) = v_{j,1} v_{j,2}$, $\pi_j(u) = u_{j,1} u_{j,2}$, $l(v_{j,1}) = l(u_{j,1})$ and either $u_{j,2}$ or $v_{j,2}$ is trivial.

Further define $f(av, au) = f(v, u)$, for any word $a$.

To verify the sufficiency of a particular set of difference histories, the following lemma is useful. We state it in some generality, since, although we only apply it to wreath products over SHORTLEX, it might also be used for wreath products over weighted orders; then we deduce the result we need for the wreath product over SHORTLEX as a corollary.

**Lemma 5.3** *Suppose that $\preceq$ is a wreath product order over* ORD, *for which any generator and its inverse have the same level. Suppose that $\mu$ is an additive positive valued function on words which is compatible with* ORD. *Let $\mathcal{D}$ be a difference machine, which is inverse closed, substring closed, and minimally labelled with respect to $\mathcal{D}$. Let $w_1, w_2$ be words which $\mathcal{D}$-fellow travel, and let $d = d_{w_1, w_2}$ be the target word difference of $(w_1, w_2)$. Then if $w_2$ is $\mathcal{D}$-irreducible, and no longer than $w_1$, and if $J$ is the maximal level of a generator in either $w_1$ or $w_2$,*

$$\mu(\pi_J(w_2)) - \mu(\pi_J(w_1)) \leq \mu(\pi_J(d))$$

PROOF: Since a generator and its inverse must have the same level, the word $w_1^{-1} w_2$ has level at most $J$, and hence, since $\mathcal{D}$ is minimally labelled, so does $d_{w_1, w_2}$. Now, provided we restrict consideration to words of level at most $J$, the composite function $\mu \circ \pi_J$ is compatible with $\preceq$. So we can apply lemma 5.1 to get the result. $\square$

**Corollary 5.4** *Suppose that $\preceq$ is a wreath product order over* SHORTLEX, *for which any generator and its inverse have the same level. Let $\mathcal{D}$ be a difference machine, which is inverse closed, substring closed, and minimally labelled with respect to $\mathcal{D}$. Let $w_1, w_2$ be words which $\mathcal{D}$-fellow travel, and let $d = d_{w_1, w_2}$ be the target word difference of $(w_1, w_2)$. Then if $w_2$ is $\mathcal{D}$-irreducible, and no longer than $w_1$, and if $J$ is the maximal level of a generator in either $w_1$ or $w_2$,*

$$l(\pi_J(w_2)) - l(\pi_J(w_1)) \leq l(\pi_J(d))$$

We are now in a position to describe and verify a finite, sufficient set of difference histories for a wreath product order over SHORTLEX.

**Proposition 5.5** *Suppose that $\preceq$ is a wreath product order over* SHORTLEX, *for which any generator and its inverse have the same level. Let $\mathcal{D}$ be a difference machine, which is inverse closed, substring closed, and minimally labelled with respect to $\mathcal{D}$. Let $K$ be the maximal number of generators of any weight in a word difference. Then the set of difference histories $H$ of the form*

$$(len, \lambda(v), \lambda(u), c_1, \ldots, c_{J-1}, c_J),$$

*as described above, such that, for each $j$, $c_j = 0, 1, -1$ or $(g(v_{j,1}, u_{j,1}), v_{j,2}, u_{j,2})$, where $v_{j,2}$ and $u_{j,2}$ both have weight at most $K$, is sufficient.*

PROOF: NB. We shall use $x <_{\text{SL}} y$ to indicate that $x$ precedes $y$ in SHORTLEX.

The proof is analogous to the proof of proposition 4.5

Suppose that $v$ is minimally $\mathcal{D}$-reducible, and that $u$ is any reduction of $v$. Without loss of generality we may assume that $v$ and $u$ have no common prefix. Let $v' = v(i)$ be a proper prefix of $v$, and $u' = u(i)$ the corresponding prefix of $u$. Let $J$ be the maximum of the levels of $v'$ and $u'$, and let

$$\delta = (len, \lambda_{v'}, \lambda_{u'}, c_1, \ldots, c_{J-1}, c_J)$$

be the difference history for $(v', u')$. Note that for any $j \le J$, if $c_j \notin \{0, 1, -1\}$, then either $\lambda_{u'} \le j$, $l(u') = l(v')$ and $\pi_j(u') <_{\text{SL}} \pi_j(v')$, or $\lambda_{v'} \le j$ and $\pi_j(v') <_{\text{SL}} \pi_j(v')$.

Suppose first that $u'$ is shorter than $v'$, and that $c_j \notin \{0, 1, -1\}$. Then $\pi_j(v') <_{\text{SL}} \pi_j(u')$, so $v'_{j,2}$ is trivial. Suppose on the other hand that $u'$ has the same length as $v'$. Then we can apply corollary 5.4 with $w_1 = u'$ and $w_2 = v'$ to see that

$$l(\pi_J(v')) - l(\pi_J(u')) \le K,$$

and hence that $v'_{J,2}$ has length at most $K$. Further, if $j < J$, where we define $m_j$ to be the largest integer such that neither $v(m_j)$ nor $u(m_j)$ contains a generator of weight greater than $j$, we can apply corollary 5.4 to the words $u(m_j)$ and $v(m_j)$ to see that

$$l(\pi_j(v(m_j))) - l(\pi_j(u(m_j))) \le K,$$

and hence that $v'_{j,2}$ has length at most $K$. (For if $v_{j,2}$ is non-trivial, it is a suffix of $\pi_j(v'(m_j))$ of length at most $l(\pi_j(v(m_j))) - l(\pi_j(u(m_j)))$. Hence we see that, for any $u$, if $c_j \notin \{0, 1, -1\}$, then $l(v'_{j,2}) \le K$.

Now suppose that $v$ is directly $\mathcal{D}$-reducible, and that $u$ is a $\mathcal{D}$-reduction of it which is $\mathcal{D}$-irreducible. Then $u'$ is $\mathcal{D}$-irreducible, and so we can apply corollary 5.4 to $u'$ and $v'$ to get

$$l(\pi_J(u')) - l(\pi_J(v')) \le K,$$

and hence deduce that if $c_J \notin \{0, 1, -1\}$ then $u'_{J,2}$ has length at most $K$. Similarly, for any $j < J$, applying corollary 5.4 to $u(m_j)$ and $v(m_j)$ we get

$$l(\pi_j(u(m_j))) - l(\pi_j(v(m_j))) \le K,$$

and hence deduce that if $c_j \notin \{0, 1, -1\}$ then $u'_{j,2}$ has length at most $K$. Hence $\delta$ is in $H$, and so $H$ carries a reduction for $v$.

So now suppose that $H$ does not carry a reduction for $v$. Then $v$ cannot be directly $\mathcal{D}$-reducible. Choose $i$ minimally so that where $v' = v(i)$ and $u' = u(i)$, the associated difference history $\delta$ is not in $H$.

For each $j < J$, a word of length $|l(\pi_j(v(m_j))) - l(\pi_j(u(m_j)))|$ is one of the components of the difference history for $(v'(m_j), u'(m_j))$. So (by the minimality of $i$), certainly $|\pi_j(v(m_j))) - l(\pi_j(u(m_j)))| \leq K$, and so $l(u_{j,2}) \leq K$ for $j < J$. Hence, since $\delta \notin H$, $u_{J,2}$ has length greater than $K$, that is

$$l(\pi_J(u')) - l(\pi_j(v'))| > K \geq l(\pi_J(d_{v',u'}))$$

Now, where $d = d_{v',u'}$,

$$l(\pi_J(v'd)) \leq l(\pi_J(v')) + l(\pi_J(d)) < l(\pi_J(u'))$$

So $v'd \prec u'$. Since $v'$ is no shorter than $u'$ and $\mathcal{D}$ is prefix closed, $v'd$ and $u'$ must $\mathcal{D}$-fellow travel. Hence $v'd$ is a direct $\mathcal{D}$-reduction of $u'$. Thus $u'$ is directly $\mathcal{D}$-reducible, and so $H$ must carry a reduction for $u'$, and also for $u$. □

## 6 Examples

### 6.1 Weighted orders

The original motivation to extend these algorithms came from Epstein's result ([6]) for geometrical finite hyperbolic groups; Lang's recent result ([14]) now shows that such groups are in fact SHORTLEX automatic, but it still seemed of interest to examine Epstein's automatic structures. Experiments described below with some hyperbolic knot groups produced automatic structures which verified Epstein's result, but the results were a little disappointing in that the structures were larger than other SHORTLEX structures.

Epstein's result suggests that a WTLEX automatic structure can be found for the fundamental group of a geometrically finite hyperbolic manifold by associating weights to the generators in such a way that low weights are assigned to generators corresponding to paths around cusps of the manifold. In a private communication, Epstein suggested that an automatic structure for a hyperbolic knot groups should therefore be found by adjoining to the generators of the Wirtinger presentation (as defined in, for example [8]) one additional generator corresponding to a path around the knot, giving that generator and one other weight one, and all the others the same, possibly higher, weight.

Three hyperbolic knots were investigated, those labelled as $4_1$, $5_2$ and $7_4$ in the table of [21] (and proved hyperbolic in [19, 20]); $4_1$ is commonly known as the figure 8 knot. With respect to the Wirtinger generators, all three were easily proved to be SHORTLEX automatic. The word acceptors had 18, 21 and 31 states, and the difference machines $D_L$ 21, 29 and 43 states respectively.

With the addition of an extra generator $a$ representing a path around the knot, presentations

$$\langle a, x, y, z, t \mid a^{-1}xy^{-1}zt^{-1},\quad x^{-1}z^{-1}tz,\quad y^{-1}txt^{-1},\quad z^{-1}x^{-1}yx \rangle$$
$$\langle a, x, y, z, t, u \mid a^{-1}tuxzy,\quad xt^{-1}x^{-1}z,\quad ty^{-1}t^{-1}x,\quad tzu^{-1}z^{-1},\quad uz^{-1}u^{-1}y \rangle$$
$$\langle a, x, y, z, t, u, v, w \mid a^{-1}vu^{-1}x^{-1}w^{-1}y^{-1}z^{-1}t^{-1},\quad t^{-1}xzx^{-1}, xvy^{-1}v^{-1},$$
$$yuy^{-1}v,\quad uyu^{-1}z^{-1},\quad wtw^{-1}u^{-1},\quad twt^{-1}x^{-1} \rangle$$

were deduced for the three knot groups from their Wirtinger presentations. In each case, the group proved in fact to be SHORTLEX automatic on the larger set of generators, with the generators ordered with $a$ first, and then $x, y, z, t \ldots$; of course this is the special case of the WTLEX order where all generators have weight 1. The automatic structures were much larger than those found using just the Wirtinger generators. The word acceptors had 781, 265 and 523 states, and the difference machines $D_L$ 287, 383 and 647 states respectively. The computation was also significantly larger; the latter two calculations in particular took 2500 seconds and 6300 seconds on a SUN Sparc 20workstation with 256 MByte RAM. (This calculation was actually done with the existing **kbmag** software for SHORTLEX, which is, of course, intrinsically more efficient than the more general algorithm for WTLEX).)

For each of the three groups, experimentation with Knuth-Bendix on WTLEX, with weight 1 assigned to $a$ and $x$ and equal weights from 2 through to 10 assigned to remaining generators, produced no evidence of stablility in the number of word differences coming from the associated rewrite system.

Certain Coxeter groups with non-standard generating sets (that is, with additional involutions introduced) were proved by Hermiller in [11] to have finite rewriting systems with respect to WTLEX. A number of examples of this type were investigated using the algorithm described in this paper and found to be WTLEX automatic, WTSHORTLEX automatic and SHORTLEX automatic, for those generating sets. In all cases investigated, the WTLEX language was the simplest of the three languages for that generating set (that is, the word acceptor had fewer states), but the SHORTLEX language for the standard generating set was simpler than that.

## 6.2  Wreath product over shortlex

In this section four infinite classes of groups are described, all of which are either automatic or asynchronously automatic (in the latter case only a more general fellow traveller property is required, see [6] for a full definition) with respect to a wreath product order over SHORTLEX. At least for small values of the parameters $p, q$, the algorithm described in this paper can be used to verify automaticity where relevant; but in each case a combinatorial proof is also straightforward to obtain.

For $p, q > 0$, the Baumslag-Solitar group $G_p, q$ is defined by the presentation

$$G_{p,q} = \langle x, y \mid yx^p = x^q y \rangle$$

Where $x, x^{-1}$ have level 1, and $y, y^{-1}$ have level 2, $G_{p,q}$ has finite rewriting system:-

$$\begin{aligned}
xx^{-1} &\to 1, & x^{-1}x &\to 1, & yy^{-1} &\to 1, & y^{-1}y &\to 1, \\
x^{r+1}y^{-1} &\to x^{r+1-p}y^{-1}x^q, & x^{r-p}y^{-1} &\to x^r y^{-1} x^{-q}, & r &= \lfloor p/2 \rfloor \\
x^{s+1}y &\to x^{s+1-q}yx^p, & x^{s-q}y &\to x^s yx^{-p}, & s &= \lfloor q/2 \rfloor.
\end{aligned}$$

The minimal representatives of group elements under that order define the regular language of words of the form

$$y^a x^{b_1} y^{c_1} \ldots x^{b_k} y^{c_k} x^d$$

where $a, k, d \in Z$, $k \geq 0$, and for each $i \leq k$,
either $c_i > 0$ and $b_i \in \{1, \ldots s, -1, \ldots, s+1-q\}$,
or $c_i < 0$ and $b_i \in \{1, \ldots r, -1, \ldots, r+1-p\}$.
This is recognised by an automaton with $\max\{p, q\} + 4$ states.

When $p = q$ the language satisfies the fellow traveller property, so we have an automatic structure. For particular values of $p, q$ this can easily be verified using the algorithm described in this paper. However, a combinatorial argument, examining the relationship between pairs of accepted words $(v, u)$ such that $u = vy$, provides a brief, alternative proof, covering all $p, q \geq 0$.

The same argument shows that when $p \neq q$, the language satisfies an asynchronous fellow traveller property, that is, for $v, w \in L$ with $w =_G vy$, and $w = zw'$, $v = zv'$, travellers on $v'$ keep close to travellers on $w'$, provided that they move at appropriate speeds (in fact, with speeds which, except at the beginning, and depending on $v'$ and $w'$, stay in the ratio $p : q$ or $q : p$). The algorithm fails to verify a synchronous fellow traveller property; in fact no such is satisfied.

In fact, in [6], Chapter 7, $G_{p,q}$ is proved to be automatic for $p = q$, and asynchronously automatic (but <u>not</u> automatic) for $p \neq q$. For $p \neq q$, $G_{p,q}$ cannot be automatic, because it has an exponential Dehn function.

Although this fact is not mentioned in [6], the language of that proof can be recognised as the set of minimal representatives of group elements under a wreath product over LEX. $G_{p,q}$ also has a finite rewriting system with respect to that order; the language of minimal representatives is very similar to the above, and the same fellow traveller properties can be similarly verified. But this order is less appropriate for computation; it admits infinite descending chains, and so is not a reduction order.

Calculation proves that $G_{p,p}$ is also SHORTLEX automatic (at least for small $p$); but wreath product orders give a much simpler language.

For $p, q > 0$, the group $G_{p,-q}$, defined by

$$G_{p,-q} = \langle x, y \mid yx^p = yx^{-q} \rangle \quad p, q, > 0$$

also has a finite rewriting system:-

$$\begin{aligned}
xx^{-1} &\to 1, \quad x^{-1}x \to 1, \quad yy^{-1} \to 1, \quad y^{-1}y \to 1, \\
x^{r+1}y^{-1} &\to x^{r+1-p}y^{-1}x^{-q}, \quad x^{r-p}y^{-1} \to x^r y^{-1} x^q, \quad r = \lfloor p/2 \rfloor \\
x^{s+1}y &\to x^{s+1-q}yx^{-p}, \quad x^{s-q}y \to x^s yx^p, \quad s = \lfloor q/2 \rfloor
\end{aligned}$$

The same language as for $G_{p,q}$ gives an automatic structure for $p = q$; the algorithm verifies this easily. A combinatorial argument provides a brief, alternative proof, that $G_{p,-p}$ is automatic for all $p > 0$ and that $G_{p,-q}$ is asynchnously automatic for $p, q > 0$ and $p \neq q$.

Again, computation shows $G_{p,-p}$ is SHORTLEX automatic at least for small $p$, but the language is much more complicated than that provided by the wreath product order.

For $p, q > 0$ the group $H_{p,q}$ defined by

$$H_{p,q} = \langle x, y \mid x^p y = y^{-1} x^q \rangle \quad p, q, > 0$$

has a finite rewriting system

$$\begin{aligned}
xx^{-1} &\to 1, \quad x^{-1}x \to 1, \quad y^{-1} \to x^p yx^{-q}, \quad yx^p y \to x^q, \\
x^{t+1}y &\to x^{t+1-(p+q)}yx^{p+q}, \quad x^{t-(p+q)}y \to x^t yx^{-(p+q)}, \quad t = \lfloor (p+q)/2 \rfloor.
\end{aligned}$$

The minimal representatives of group elements under that order define the regular language of words of the form

$$x^{a_1}y^{b_1} \ldots x^{a_k}y^{b_k}x^c$$

where $a, k, d \in Z$, $k \geq 0$, for each $i \leq k$, $b_i > 0$, $a_1 < p$, and for $1 < i \leq k$,

$$a_i \in \{1, \ldots, t, -1, \ldots, -(p+q) + t + 1\}$$

The fellow traveller property holds for all $p, q$; this can easily be verified combinatorially, and essentially follows from the relation

$$x^{p+q}y = yx^{p+q}$$

(which shows that $H_{p,q}$ is a quotient of $G_{p+q,p+q}$.) So $H_{p,q}$ is automatic for all $p, q$.

For $p, q > 0$ the group $H_{p,-q}$ defined by the presentation

$$H_{p,-q} = \langle x, y \mid x^p y = y^{-1} x^{-q} \rangle \quad p, q, > 0$$

has finite rewriting system

$$xx^{-1} \to 1, \quad x^{-1}x \to 1, \quad y^{-1} \to x^p y x^q, \quad y x^p y \to x^{-q}.$$

The minimal representatives of that order define the regular language of words of the form
$$x^{a_1} y^{b_1} \ldots x^{a_k} y^{b_k} x^c$$
where $a, k, d \in Z$, $k \geq 0$
for each $i \leq k$, $b_i > 0$, $a_1 < p$, and for $1 < i \leq k$, $a_i \neq p$.

The fellow traveller property holds for all $p, q$. This can very easily be verified combinatorially; if $u, v \in L$, and $l(v^{-1}u) \leq 1$, then the words $u, v$ must be of the form $wu_0, wv_0$ for short words $v_0, u_0$. So in fact, the language is a *broomlike combing*, in the sense of [1]; hence by [1] this group is virtually free.